\documentclass[12pt,reqno,a4paper]{amsart}

\usepackage{tikz}
\usetikzlibrary{arrows,backgrounds,calc}
\usepackage{calrsfs}
\usepackage[charter]{mathdesign}
\usepackage{amsthm, amsmath, amsfonts}
\usepackage{graphicx}
\usepackage{hyperref}
\usepackage{lineno}
\usepackage{enumerate}
\usepackage{mathtools}
\usepackage{pgfplots}
\pgfplotsset{compat=1.5}

\title{A bijection between phylogenetic trees and plane oriented recursive trees}

\author[H.~Prodinger]{Helmut Prodinger}
\address[Helmut Prodinger]{Department of Mathematical
	Sciences, Stellenbosch University, 7602 Stellenbosch,
	South Africa}
\email{\href{mailto:hproding@sun.ac.za}{hproding@sun.ac.za}}

\begin{document} 
\begin{abstract}
	 Phylogenetic trees are binary nonplanar trees with labelled leaves, and plane oriented recursive trees are planar trees with an increasing labelling. Both families are enumerated by double factorials. A bijection is constructed, using the respective representations a 2-partitions and trapezoidal words. 
\end{abstract}

\maketitle

\section{Introduction}

There are many occurrences of the \emph{double factorials}
 $$(2n-1)!! =1\cdot3\cdots(2n-1)$$
  in the combinatorial literature. A nice survey is by Callan~\cite{Callan09}.

Two manifestations of them deal with trees, and it is our objective to establish a bijection between them.

The phylogenetic trees are binary nonplanar trees with the leaves labelled by the numbers $1,2,\dots,n+1$.
Their number is given by $(2n-1)!!$. Stanley describes (codes) them by set partitions of $\{1,2,\dots,2n\}$
into $n$ sets of $2$ elements each. These are easily enumerated by the double factorials: Just note that they are counted by $\frac{(2n)!}{2^nn!}$, where we start with all permutations of $2n$ elements, but divide by $n!$, since the order of the blocks does not count, and by $2^n$, since in each block the order of the 2 elements is immaterial. We call these set partitions 2-partitions.  Stanley's coding will be reviewed in a later section.

Plane oriented recursive trees (PORTs)~\cite{Hwang}, also known as heap ordered trees, are planar trees, where the nodes are labelled by the integers $1,\dots,n+1$, and the labels are increasing towards the leaves. They are also enumerated by  $(2n-1)!!$. They are also coded by simple objects called \emph{trapezoidal words}, which is reviewed in the next section. 

Our bijection will in fact be between the two codings, i.e., between 2-partitions and trapezoidal words.

\section{Trapezoidal words and Plane oriented recursive trees}

One of the easiest manifestations of double factorial is by words\\ $x_1x_2\dots x_n$, where $1\le x_i\le 2i-1$;
they were called \emph{trapezoidal words} by Riordan~\cite{Riordan76}; see also \cite{Callan09}.

Plane oriented recursive trees (PORTs) are rooted planar trees, where the $n$ nodes are labelled by the numbers $1,\dots,n$ in an increasing way from the root to the leaves.

If one has already  such a PORT with $n$ nodes, there are $2n-1$ positions where a new node labelled $n+1$ can be attached, whence the enumeration by double factorials: The number of PORTs with $n+1$ nodes is given by $(2n-1)!!$, and the trees are in obvious bijection with trapezoidal words of length $n$, simply by recording the position where one node after the other was inserted. PORTs were also called heap ordered trees, but we adopted the notation from \cite{Hwang}.

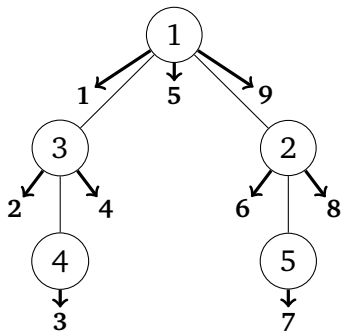
\begin{figure}[h]\centering
\begin{tikzpicture}[level distance=1.5cm,
level 1/.style={sibling distance=1.5cm},
level 2/.style={sibling distance=1.5cm}]
\tikzstyle{every node}=[circle,draw]
\tikzstyle{solid node}=[circle,draw,inner sep=1.5,fill=black]
\tikzstyle{pointed node}=[level distance=0.8cm, sibling distance=0.6cm, edge from
parent/.style={draw,very  thick, ->},
every child node/.style={inner sep=0.5pt, draw=none, font=\footnotesize\bf}]

\node (Root) {1}
child[pointed node] {node {1}}
child { node {3}  
  child[pointed node] {node {2}}
  child { node {4}
  	  {child[pointed node] {node {3}}} }
  child[pointed node] {node {4}}
}
child[pointed node] {node {5}}
child { node {2} {child[pointed node] {node {6}}child {
		 node {5} {child[pointed node] {node {7}}}
		 }child[pointed node] {node {8}}}}
	child[pointed node] {node {9}};
\end{tikzpicture}
\caption{A PORT with 5 nodes and the 9 positions where node labelled 6 could be inserted.}
\end{figure}

\section{2-partions and Phylogenetic trees}

Phylogenetic trees are nonplanar binary trees, with the leaves labelled by the numbers $1,\dots,n+1$ in an arbitrary way. Stanley \cite{Stanley99} describes the  procedure in Figure~2 to label the remaining nodes as well (except for the root): The numbers $n+2,\dots,2n$ are used as labels in this order as follows: among the pairs of children that are both labelled, but the parent isn't, find the smallest label of a child; it is the parent who gets the current label. The procedure can be seen from   Figure~2. At the end, the labels of each pair of 2 children form a subset of 2 elements, leading to $n$ such pairs and thus to a 2-partition.
See also \cite{Bona} for more results and references about phylogenetic trees.

\begin{figure}[h]\label{F1}
	\begin{tikzpicture}[scale=0.8,level distance=1.5cm,
	level 1/.style={sibling distance=2.5cm},
	level 2/.style={sibling distance=1.5cm}]
	\tikzstyle{every node}=[circle,draw,font=\scriptsize,inner sep=1.5pt]

	\node (Root) {\phantom{1}}
	child {node {\phantom{1}}
	    child{node{\phantom{1}}
	    	child{node{1}}	
	    	child{node{4}}}	
		child{node{6}}
		}
	child {node {\phantom{1}}
		child{node{\phantom{1}}child{node{2}}	
			child{node{\phantom{1}}child{node{5}}	
				child{node{7}}
				}}	
		child{node{3}}}
	;
	\end{tikzpicture}\hspace*{1cm}
\begin{tikzpicture}[scale=0.8,level distance=1.5cm,
level 1/.style={sibling distance=2.5cm},
level 2/.style={sibling distance=1.5cm}]
\tikzstyle{every node}=[inner sep=2.3pt,circle,draw,font=\scriptsize]
\tikzstyle{dobo}=[inner sep=1.0pt,circle,draw,font=\scriptsize]

\node (Root) {\phantom{1}}
child {node[dobo]{12}
	child{node{8}
		child{node{1}}	
		child{node{4}}}	
	child{node{6}}
}
child {node[dobo] {11}
	child{node[dobo]{10}child{node{2}}	
		child{node{9}child{node{5}}	
			child{node{7}}
		}}	
		child{node{3}}}
	;
	\end{tikzpicture}
			\caption{Left: A nonplanar binary tree with leaves labelles by $1,\dots,7$.
							Right: the remaining nodes (except for the root) are now labelled by $8,\dots,12$.  The two children of each node form the 2-partition: $\{1,4\}, \{2,9	\}, \{5,7\}, \{6,8\}, \{3,10\}, \{11,12\}$}
\end{figure}
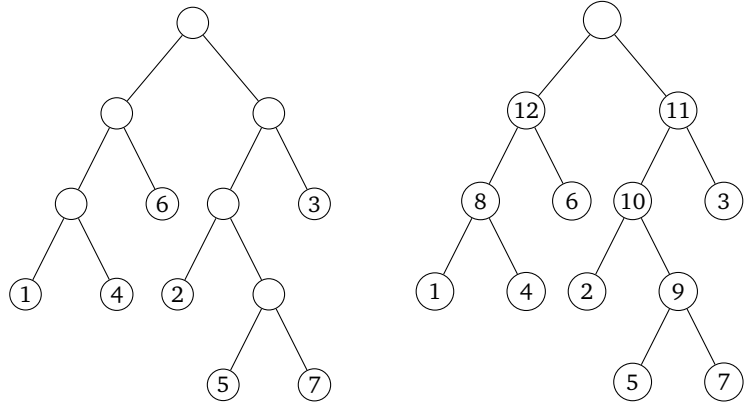

Although Stanley leaves it to the reader to figure out why this works, we sketch a possible answer by showing how a phylogenetic tree can be reconstructed from a 2-partition:

We use the consecutive labels $n+1,\dots,2n$ to work as a parent. For that, we look at the block, such that both entries are smaller than the current new label, and among them at the one in which the smallest number occurs.
After that,  the pair is discarded, and the process continues until all pairs have been processed. The final root may be thought of having the label $\infty$. 

Thus, in the running example, the number 8 is the current new parent, and $\{1,4\},  \{5,7\}$ are such that both members are smaller than 8. The block $\{1,4\}$ is chosen. Then we move to number $9$.  Candidates are $\{5,7\}, \{6,8\}$; the block containing the number 5 is used, then discarded, and so on.

\section{The bijection between 2-partitions and trapezoidal words}

Our strategy of proof is to grow a 2-partition of $2n-2$ elements to one of $2n$ elements (hereby establishing once again the $(2n-1)!!$ formula), and describing how the corresponding trapezoidal word of length $n-1$ grows to one of length $n$. The correspondence is bijective at each step. Our argument is essentially by induction. It should be noted that the way a 2-partition (and a trazoidal word) grows towards a final partition/word is \emph{unique}. 

Two new elements $2n-1,2n$ can form a class of their own, and this can happen in $(2n-3)!!$ ways. Otherwise,
 $2n$ matches with an element $1\le b\le 2n-2$ (in $2n-2$ ways), and $2n-1$ matches with the former partner $a$ of $b$. So, the set $\{a,b\}$ is replaced by the sets $\{a,2n-1\}, \{b,2n\}$. Such an operation is often called a \emph{rotation}. Thus, we get altogether $(1+(2n-2))\cdot (2n-3)!!=(2n-1)!!$ 2-partitions, as to be expected.

Then we augment the corresponding trapezoidal word $x_1\dots x_{n-1}$, by $x_n$, defined a follows:
if the second case happened and $2n$ matches with an element $1\le b\le 2n-2$, then we set $x_n:=b$, otherwise, if $2n-1,2n$  form a class of their own we set $x_n:=2n-1$. It is easy to see that this operation is reversible.

As an example, let us see how the trapezoidal word $1,2,5,5,2,4$ lets the 2-partition grow. The first 1 translates into the starting partition $\{1,2\}$.
\begin{gather*}
\{1,2\}\xrightarrow{2}\{1,3\}\{2,4\}
\xrightarrow{5}
\{1,3\}\{2,4\}\{5,6\}
\xrightarrow{5}
\{1,3\}\{2,4\}\{6,7\}\{5,8\}\\
\xrightarrow{2}
\{1,3\}\{2,10\}\{4,9\}\{6,7\}\{5,8\} \\
\xrightarrow{4}
\{1,3\}\{2,10\}\{4,12\}\{9,11\}\{6,7\}\{5,8\} 
\end{gather*}

And here is how the PORT develops:
\begin{figure}[h] 
	\begin{tikzpicture}[level distance=1.5cm,
	level 1/.style={sibling distance=1.5cm},
	level 2/.style={sibling distance=1.5cm}]
	\tikzstyle{every node}=[circle,draw,font=\footnotesize]
	
	\hspace*{-4cm};
	\node (Root) {1};\hspace*{1cm}
\node (Root) {1}
child { node {2} }
;
\hspace*{1cm}
\node (Root) {1}
child { node {2}child { node {3} } }
;
\hspace*{1.8cm}
\node (Root) {1}
child { node {2}child { node {3} } }
child{node{4}}
;
\hspace*{3.2cm}
\node (Root) {1}
child { node {2}child { node {3} } }
child{node{5}}
child{node{4}}
;
	\end{tikzpicture}

\end{figure}

\begin{figure}[h] 
	\begin{tikzpicture}[level distance=1.5cm,
	level 1/.style={sibling distance=1.5cm},
	level 2/.style={sibling distance=1.5cm}]
	\tikzstyle{every node}=[circle,draw,font=\footnotesize ]
	
	\hspace*{-2.5cm};
	\node (Root) {1}
	child { node {2}child { node {6} }child { node {3} } }
	child{node{5}}
	child{node{4}}
	;
	\hspace*{5cm};
	\node (Root) {1}
	child { node {2}child { node {6} }child { node {7} }child { node {3} } }
	child{node{5}}
	child{node{4}}
	;
	\end{tikzpicture}
	
\end{figure}

\bibliographystyle{plain}

\end{document}